\documentclass[10pt,A4,conference]{IEEEtran}
\usepackage{geometry}
 \usepackage{authblk}
\usepackage{tikz}  
\usetikzlibrary{calc,intersections}
 \usepackage{verbatim}
\usepackage{textcomp}
\usepackage[T1]{fontenc}
\usepackage{bbm}
\usepackage{multirow}
\usepackage{dsfont}
\usepackage{cite}
\usepackage{epsfig}
\usepackage{float}
\usepackage{enumerate}
\usepackage{balance}
\usepackage{color}
\usepackage{subcaption}
\usepackage{caption}
\usepackage{algorithm}
\usepackage[noend]{algpseudocode}
\algrenewcommand\Return{\State \algorithmicreturn{} } 
\usepackage{hyperref}
\hypersetup{
    colorlinks=true,
    linkcolor=blue,
    filecolor=magenta,      
    urlcolor=cyan,
		citecolor=blue,
}
\usepackage{amssymb}
\usepackage{amsthm}
\usepackage{array}
\usepackage{graphicx}
\usepackage{lettrine}
\usepackage{epsf} 
\usepackage{psfrag}
\usepackage[usenames,dvipsnames]{pstricks}
\usepackage{pst-grad} % For gradients
\usepackage{pst-plot} % For axes
\usepackage[space]{grffile} % For spaces in paths
\usepackage{etoolbox} % For spaces in paths
\usepackage{url} 

\usepackage{multirow}

\usepackage{optidef}
\usepackage[capitalize]{cleveref}
\usepackage[colorinlistoftodos,bordercolor=orange,backgroundcolor=orange!20,linecolor=orange,textsize=scriptsize]{todonotes}
\Crefname{figure}{\Fig.}{Figs.}

\begin{document}

\title{Muti-Agent Proximal Policy Optimization For Data Freshness in UAV-assisted Networks}

 %\author{\authorblockN{Mouhamed Naby Ndiaye\authorrefmark{1}, Hajar El Hammouti\authorrefmark{1}, ElHoucine Bergou\authorrefmark{1} }
 %\authorblockA{\authorrefmark{1}School of Computer Science, Mohammed VI Polytechnic University (UM6P), Benguerir, Morocco}}
 
\author[1]{Mouhamed Naby Ndiaye}
\author[1]{El Houcine Bergou}
% \author[2]{Mounir Ghogho}
\author[1]{Hajar El Hammouti}
\affil[1]{School of Computer Science, Mohammed VI Polytechnic University (UM6P), Benguerir, Morocco, \authorcr
emails: {\{naby.ndiaye,elhoucine.bergou,hajar.elhammouti\}@um6p.ma}}
% \affil[2]{TICLab, International University of Rabat (UIR), Rabat, Morocco,\authorcr
% email: {mounir.ghogho@uir.ac.ma}}

%  El Houcine Bergou\authorrefmark{1}, Mounir Ghogho \authorrefmark{2}, Hajar El Hammouti\authorrefmark{1} }
%  \authorblockA{\authorrefmark{1}School of Computer Science, Mohammed VI Polytechnic University (UM6P), Benguerir, Morocco,\\\authorrefmark{2}TICLab, International University of Rabat (UIR), Rabat, Morocco,}
% emails: \{naby.ndiaye,hajar.elhammouti,elhoucine.bergou\}@um6p.ma, mounir.ghogho@uir.ac.ma}

\maketitle

% ===================
% # I. Abstract #
% ===================
\begin{abstract}
Unmanned aerial vehicles (UAVs) are seen as a promising technology to perform a wide range of tasks in wireless communication networks. In this work, we consider the deployment of a group of UAVs to collect the data generated by IoT devices. Specifically, we focus on the case where the collected data is time-sensitive, and it is critical to maintain its timeliness. Our objective is to optimally design the UAVs’ trajectories and the subsets of visited IoT devices such as the global Age-of-Updates (AoU) is minimized. To this end, we formulate the studied problem as a mixed-integer nonlinear programming (MINLP) under time and quality of service constraints. To efficiently solve the resulting optimization problem, we investigate the cooperative Multi-Agent Reinforcement Learning (MARL) framework and propose an RL approach based on the popular on-policy Reinforcement Learning (RL)  algorithm: Policy Proximal Optimization (PPO). Our approach leverages the centralized training decentralized execution (CTDE) framework where the UAVs learn their optimal policies while training a centralized value function. Our simulation results show that the proposed MAPPO approach reduces the global AoU by at least a factor of $1/2$ compared to conventional off-policy reinforcement learning approaches.

%We formulate our problem as a multi-objective optimization problem where the main objectives are to maximize the total number and minimize the Age of Information (AoI) of served IoT devices. In order to efficiently solve our optimization problem, we propose a cooperative Multi-Agent Reinforcement Learning (MARL) algorithm based on the popular on-policy RL algorithm Policy Proximal Optimization (PPO). Simulation results show that our proposed algorithm performs well compared to off-policy MARL algorithm.
\end{abstract}
\IEEEoverridecommandlockouts
\IEEEpeerreviewmaketitle

% ===================
% # I. Introduction #
% ===================

\section{Introduction}

Unmanned aerial vehicles (UAVs), also known as drones, have been considered as a promising technology to enhance network coverage and provide on-demand connectivity. While UAVs present many advantages, their deployment in a multi-agent environment still faces a number of challenges including resource allocation, trajectory design and cooperation between the drones~\cite{bajracharya20226g}. In this paper, we are interested in the case where UAVs are used to periodically collect data from Internet of Things (IoT) devices and send it to a server for analysis and decision-making. Specifically, we study the scenario where the collected data is time-sensitive, and it is critical to maintain its timeliness. To measure the freshness of the data, we adopt the \textit{age-of-updates} (AoU) metric, which captures the time since the last data update was collected. The main objective of deploying UAVs, in this context, is to minimize the global AoU of the network given a period of time. Since the UAVs operate in a common environment, it is extremely relevant to study their behaviors and interactions to achieve a common objective (i.e., minimizing the global AoU) while they perform their tasks according to their own interests and knowledge.

%To this end, a number of UAVs is deployed and their trajectories need to be designed such that the total collected data is kept as fresh as possible. In this context, we answer the question: how UAVs can efficiently interact in a multi-drone system such that the total AoU is minimized? 

%Our objective is to design the UAVs trajectory such that the total collected data is kept as fresh as possible.
%To measure the freshness of data, we adopt a new metric that has been introduced recently, namely \textit{the age-of-updates} (AoU).

%Unmanned aerial vehicles (UAVs) offer a wide range of possibilities in the field of telecommunications due to their flexibility and maneuverability~\cite{bajracharya20226g}. In fact, they can be used for multiple applications, including enhancing network coverage and collecting data from ground Internet of Things (IoT) devices. In this paper, we are interested in the case where UAVs periodically collect data from IoT devices and send it to a server for analysis and decision-making. For scenarios where the data is time-sensitive, it is important to carefully design the UAVs' trajectories so that maximum data is collected, but also  the age of the collected updates is minimized. In this work, we present a practical framework for deploying autonomous UAVs that can adapt to unknown environments and work cooperatively to minimize the age of the collected data captured by a metric called Age-of-Update (AoU). 

\subsection{Related Work}

Early studies of UAV-assisted networks assumed a single UAV setup where a drone is deployed to achieve a given task, such as improving the network coverage or maximizing the amount of collected data~\cite{mozaffari2016unmanned,mozaffari2017wireless}. In general, this class of work models the UAV task as a mixed-integer nonlinear programming (MINLP) which is solved using traditional MINLP optimization approaches (e.g., Branch and Bound and dynamic programming) or some  heuristics~\cite{bushnaq2019aeronautical,Wang2022JointOO}. For example, in~\cite{bushnaq2019aeronautical}, the authors study the trajectory of a drone such that the traveling and hovering times are minimized while the number of hovering regions is maximized. They formulate their optimization as a standard travelling salesman problem and use a heuristic approach to balance the tradeoff between the traveling and hovering times. In~\cite{Wang2022JointOO}, the authors jointly optimize the UAV trajectory and sensors' energy under power and time constraints. The studied problem is decomposed into subproblems which are solved approximately using convex optimization. While such approaches provide efficient solutions, they scale poorly with high dimensional spaces, and therefore, they are not adapted to multi-agent systems where the optimization space is intrinsically large. %Second, they are not well adapted to multi-agent systems.  

To overcome this problem, machine learning (ML) approaches, in particular, reinforcement learning (RL) based methods, have been proposed recently. The unprecedented popularity of RL stems from its ability to solve complex sequential decision-making problems where the agents interact with the environment to learn an optimal policy. Specifically, in the context of UAVs, deep RL (DRL) was proposed as a powerful tool that combines deep neural networks (DNNs) and RL to model continuous action spaces and their unlimited corresponding rewards. For example, in~\cite{Liu2018EnergyEfficientUC}, the authors minimize the energy of a drone using a DRL under time and connectivity constraints. The proposed DRL approach uses a DNN to approximate the Q-function. Similarly, in~\cite{Xiong2021UAVAssistedWE}, the authors use DRL to model the large state space of a UAV-assisted communication where the drone is charged using energy transfer. Despite its effectiveness in dealing with high-dimensional spaces, DRL is not suitable for multi-UAVs systems where the policies of the drones are interdependent. In fact, in a single-agent setting, the environment is in general stationary and the UAV can fully observe its actions and rewards. Conversely, in a multi-agent system, the drones operate simultaneously, hence, the environment experienced by each drone becomes non-stationary.

In this context, the multi-agent RL framework (MARL) claims to address these challenges by assuming that the agents make decisions based on local observations only. In this perspective, a handful of works have been proposed on the cooperation of UAVs to accomplish a common mission, usually based on MARL~\cite{Cui2020MultiAgentRL, Seid2021MultiAgentDF,9900429}. In~\cite{Seid2021MultiAgentDF}, the authors study the deployment of a multi-drone system to achieve computational tasks in an IoT network. The problem is formulated as a stochastic game, with the goal of minimizing the long-term computational cost, and is solved using a multi-agent deep reinforcement learning (MADRL) method. In~\cite{9900429}, the authors propose a cooperative multi-agent algorithm with an actor-critic architecture where UAVs and ground base stations cooperate to maximize the global cache hit ratio. They also show that the proposed approach overcomes a standard RL method. 

Most existing MARL solutions for multi-UAV systems assume value-based RL algorithms that use value functions to determine the expected reward for an action in a given state. In value-based RL, a value function is estimated for each state and the expected reward is computed as the discounted sum of future rewards. For instance, in~\cite{hu2021distributed}, the authors propose a value decomposition network (VDN) algorithm to optimize the UAVs' trajectories such that the number of served users is maximized. The proposed approach involves distributed multi-agent training where the drones learn their optimal strategies in a fully distributed fashion without having to share their local rewards. Unlike existing works, we leverage a policy-based MARL algorithm, which is built on the centralized training and decentralized execution (CTDE) paradigm where the UAVs learn their decisions based on local observations as well as some globally shared information. Policy-based RL methods estimate the value of a policy while using it for control. They also have some advantages over value-based algorithms including the ability to handle continuous action spaces and achieve faster convergence. Finally, our proposed MARL approach involves a multi-agent actor-critic architecture where the actor (a DNN that predicts optimal actions) is decentralized whereas the critic (a DNN that approximates the value function) is centralized. The training of the actor-critic architecture is performed using the policy proximal optimization (PPO) algorithm; a powerful policy gradient RL algorithm released by Open AI to handle large models and support parallel implementations~\cite{schulman2017proximal}. 

\subsection{Contribution}

We consider a UAV-assisted network where a team of drones is deployed to collect data updates from a set of IoT devices. Specifically, we study the case where the collected data is time-sensitive, and it is critical to maintain its freshness. Our objective is to optimally  design the UAVs' trajectories and the subsets of
visited IoT devices such as the global AoU is minimized. Our contributions can be summarized as follows.
\begin{itemize}

  %  \item  We consider a practical UAV-assisted network system in which a team of UAVs is deployed to act as a relay between low-capacity IoT devices and remote BS in order to collect as fresh data generated by these IoT devices on the ground as possible, despite strict mission time constraints and limited information on their environments. The UAVs can collect data generated by UEs while adapting their trajectory to minimize the AoU of data collected by all UAVs. 
%\item We propose a model of IoT network where devices generate data with heterogeneous data update frequencies.
\item We introduce an AoU metric that takes into consideration the frequency at which the data is generated by each device. The proposed AoU also involves the UAVs trajectories together with the subsets of visited devices.
\item We formulate the studied problem as a MINLP where the objective is to optimize the global AoU under time and quality of service constraints. To solve this challenging problem, we leverage the MARL framework and propose a multi-agent PPO (MAPPO) based approach, referred to as MAPPO-AoU, where UAVs learn a centralized value function while searching their optimal policies. 
\item Our simulation results show that MAPPO-AoU requires fewer iterations to achieve convergence compared to conventional Value-based RL algorithms. Furthermore, during the execution, the proposed approach reduces the global AoU by a factor of $1/2$ compared to Value-based RL.

\end{itemize}

\subsection{Organization}
The remainder of the paper is organized as follows. In Section~\ref{Sys}, we describe the studied system model and introduce the AoU metric. The mathematical formulation of the problem is given in Section~\ref{Prob}. In Section~\ref{algo}, we describe in details the proposed MARL approach. Next, in Section~\ref{simu}, we assess the performance of the proposed approach. Finally, concluding remarks are provided in Section~\ref{Conc}.

% =============================================
% # II. System model #
% =============================================

\section{System Model}\label{Sys}
Consider a set of $I$ IoT devices that generate data updates periodically during a time span $T$. Let $\mathcal{I}$ be the set of devices. We divide $T$ into equal intervals of length $\tau$. Each device $i$ generates new data with a period $k_i\tau$, where $k_i\in\{1,\dots,K\}$, and $K=\frac{T}{\tau}$. The generated data is cumulated in a buffer until it is collected by a UAV. We assume that devices' buffers are large enough to save all the generated data during the time span $T$. Fig.~\ref{Fig1} illustrates an example of the studied system model. 

% We denote by $\mathcal{U}$ the set of $U$ UAVs that are deployed to collect data updates from ground devices.  and transmitted to a distant server
\subsection{Age-of-Updates Metric}
\begin{figure}
    \centering
   \includegraphics[width=\linewidth]{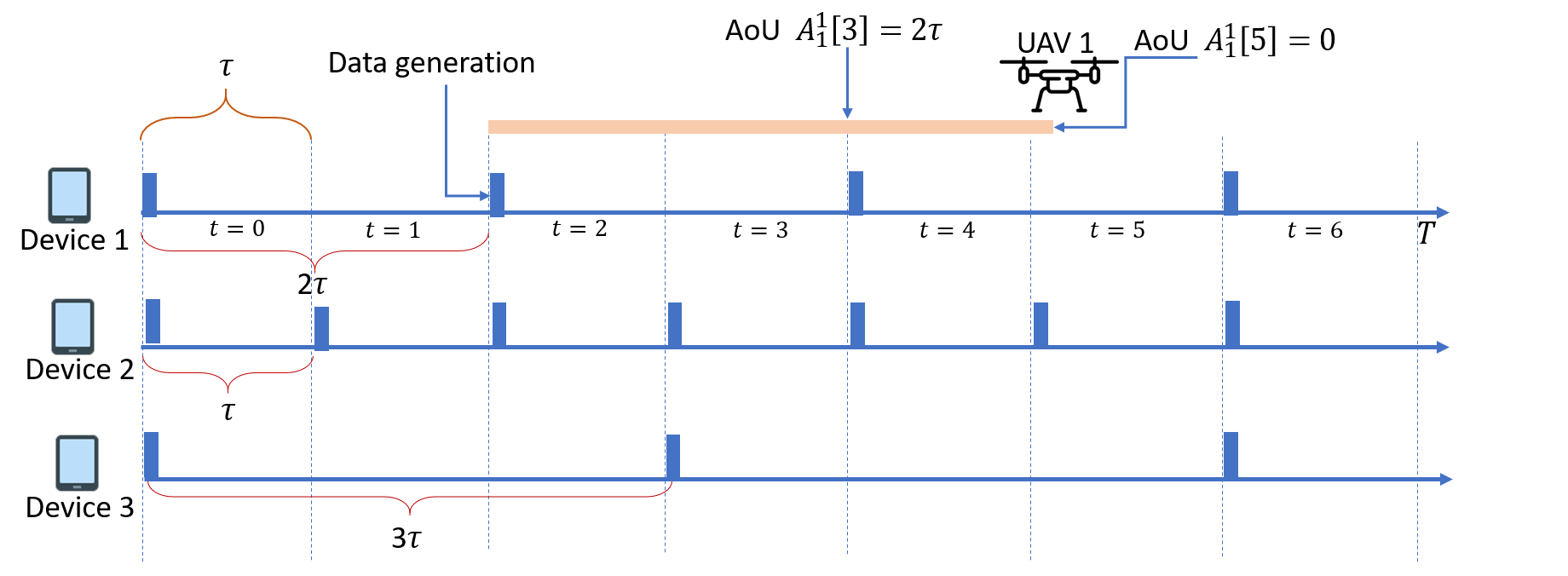}
    \caption{An example of $3$ devices with different data generation periods. Devices~$1$, $2$ \\ and $3$ generate data updates each $2\tau$, $\tau$ and $3\tau$, respectively. In this example, the time span $T$ is divided into $K=7$ time intervals of length $\tau$. An example of AoU notation is\\ given in the figure. $A^1_1[3]$ denotes the AoU of the data, generated at the period number $1$ of device $1$, during time interval $3$. In this example, $A^1_1[3]$ is equal to $2\tau$ in this\\ example. When UAV~$1$ collects the data at time interval~$5$, AoU becomes $A^1_1[5]=0$. }
    \label{Fig1}
\end{figure}

 To capture the freshness of the data, we consider the age-of-updates (AoU) metric. AoU measures the elapsed time  between the data generation and data collection. We denote by $A_{i}^n[t]$ the AoU during time interval $t \in \{1,\dots,K\}$ of the data generated by IoT device $i$ at its $n^{\rm th}$ period (i.e., during $nk_i^{\rm th}$ time interval). $A_{i}^n[t]$ can be expressed in a recursive way, as follows

\begin{equation}\label{equa1}
    A_{i}^n[t]=\Bigg\{
    \begin{array}{ll}
   A_{i}^{n}[t-1]+\tau & \text{if} \sum\limits_{l=nk_i}^{t}\sum\limits_{u \in \mathcal{U}}\alpha_{iu}[t]=0 \text{ \& } nk_i\leq t\\
    0 & \text{otherwise},
    \end{array}
\end{equation}
where $\mathcal{U}$ is the set of UAVs, $\alpha_{iu}[t]$ is a binary variable, with $\alpha_{iu}[t]=1$ if UAV $u$ collects data of device $i$ during time interval $t$ and $\alpha_{iu}[t]=0$ otherwise.  Specifically, when updates from device $i$ are not collected during time interval $t$ (i.e., $\alpha_{iu}[t] = 0$), the AoU is increased by $\tau$. Conversely, when updates are collected (i.e., $\alpha_{iu}[t]=1$) or not generated yet  (i.e., $nk_i\geq t$), the AoU is set to zero. This definition of AoU is more generic than the definition of information age given in\cite{10001424}. Indeed, these definitions of data freshness can be derived from our definition by taking $\tau=1$ and not considering any data generation period. An example of AoU notation is given in Fig.~\ref{Fig1}.
% \cite{yang2020age,AbdElmagid2019DeepRL,ourwork}
% An alternative way to write equation~(\ref{equa1}) is
% \begin{equation}
%  A_{i}^n[t]=\Big(A_{i}^n[t-1]+\tau\Big)\Big(1-\sum \limits_{u \in \mathcal{U}}\alpha_{iu}[t]\Big)\mathbbm{1}_{nk_i\leq t},
% \end{equation}
% where $\mathbbm{1}_{nk_i\leq t}$ is equal to $1$ if $nk_i\leq t$ and $0$ otherwise. 
% An example of AoU notation is given in Fig.~\ref{Fig1}. 
During $T$, the AoU of all the data generated by IoT device $i$, $f_i(\boldsymbol{\alpha}_i)$, is given by 

%Our aim is to minimize the total AoU of IoT devices during the time span $T$ given by

% \begin{equation}
% f_i(\boldsymbol{\alpha}_i)= \sum 
%    \limits_{t=1}^K \sum\limits_{n=0}^{\lfloor \frac{K}{k_i} \rfloor}\Big(A_{i}^n[t-1]+\tau\Big)\Big(1-\sum \limits_{u \in \mathcal{U}}\alpha_{iu}[t]\Big)\mathbbm{1}_{nk_i\leq t},
% \end{equation}
\begin{equation}
f_i(\boldsymbol{\alpha}_i)= \sum 
   \limits_{t=1}^K \sum\limits_{n=0}^{\lfloor \frac{K}{k_i} \rfloor}A_{i}^n[t],
\end{equation}
with $\boldsymbol{\alpha}_i=(\alpha_{iu}[t])_{\substack{u\in \mathcal{U} \\ t \in \{0,\dots,K-1\}}}$ a $U\times K$ matrix.

%\begin{equation}
%f_1(\boldsymbol{\alpha})=\sum \limits_{i \in \mathcal{I}}\sum \limits_{u \in \mathcal{U}} \sum 
%   \limits_{t=1}^K \sum\limits_{n=0}^{\lfloor \frac{K}{k_i} \rfloor}\Big(A_{iu}^n[t-1]+\tau\Big)\Big(1-\alpha_{iu}[t]\Big)\mathbbm{1}_{nk_i<t}
%\end{equation}
%with $\mathbbm{1}_{nk_i<t}$ is equal to $1$ if $nk_i<t$ and $0$ otherwise.

\subsection{Communication Model}
To send their data updates, an IoT device communicates with a UAV via the air-to-ground channel. We assume a block Rician-fading model where the propagation channels between any device-UAV pair are constant during a time interval of length less than or equal to $\tau$. Hence, the air-to-ground channel between device $i$ and UAV $u$ is modeled as follows  

\begin{equation*}
    h_{iu}[t]=d_{iu}^{-1}[t]\left(\sqrt{\frac{\Phi}{\Phi+1}} {\xi}^{\rm LoS}_{iu}[t]+\sqrt{\frac{1}{\Phi+1}} {\xi^{\rm NLoS}}_{iu}[t]\right),
\end{equation*}
where $\Phi$ is the Rician factor, ${\xi}^{LoS}_{iu}[t]$ is the line-of-sight (LoS) component with $\left|{\xi}^{\rm LoS}_{iu}[t]\right|=1$, and ${\xi}^{\rm NLoS}_{iu}[t]$ the random non-line-of-sight (NLoS) component, which follows a Rayleigh distribution with mean zero and variance one. Finally,  $d_{iu}[t]$ is the distance between the pair device-UAV during time interval $t$ that is given by 
\begin{equation*}
    d_{i,u}[t]=\sqrt{(x_{u}[t]-x_{i})^2+(y_{u}[t]-y_{i})^2+(H_{u})^2},
\end{equation*}
 where $(x_u[t],y_u[t],H_u)$ and $(x_i,y_i,0)$ are the 3D positions of UAV $u$ and device $i$ respectively.  To avoid collisions, each UAV $u$ flies at a constant and different altitude  $H_u$.

Assume devices use orthogonal frequency division multiple access (OFDMA) for their
transmissions. Accordingly, the signal-to-noise ratio (SNR) of IoT device $i$ with respect to UAV $u$ is given by
 \begin{equation*}
     \Gamma_{iu}[t]= P_{i}[t]\left|h_{iu}[t]\right|^{2} / \sigma^{2}, 
 \end{equation*}
 where $P_i[t]$ is the transmit power of device $i$ during time interval $t$, and  $\sigma^{2}$ is the variance of an additive white Gaussian noise. Therefore, the rate of IoT device $i$ with respect to UAV $u$ during time slot $t$ can be written
 
 \begin{equation*}
     {R}_{iu}[t]=B_{iu}[t]\log _{2}\left(1+\Gamma_{iu}[t]\right),
 \end{equation*}
where $B_{iu}[t]$ is the allocated bandwidth between device $i$ and UAV $u$ during time slot $t$. 

 Each UAV starts its flight from its charging and docking station (CDS) located at $(x_{u}^{\rm CDS},y_{u}^{\rm CDS},0)$ and goes back to the CDS before the depletion of its battery. During its flight, a UAV stops to collect data from subsets of IoT devices. To ensure a successful and rapid data transmission between device $i$ and UAV $u$, the data rate between the pair device-UAV, $R_{iu}$, should be greater than a predefined threshold $R^{\rm min}$, where $R^{\rm min}$ is high enough to ensure a quasi-instant transmission of all the accumulated data updates\footnote{This is also possible since the size of the updates for some applications (e.g., environmental sensors) is small and can be rapidly transmitted to the UAV.}. Accordingly, the data collection time is negligible compared to the flight time\footnote{The case where the collection time is not negligible  is left as part of our future work.}. We also assume that  UAVs fly at a constant speed $V$. They first fly vertically to reach the target altitude, collect data while flying in the horizontal plane, and then land vertically on the CDS.
 Thus, the total flight time of a UAV $u$, $\zeta_u$, is given by 

\begin{align}
 &\zeta_{u}(\boldsymbol{x}_u,\boldsymbol{y}_u)=2\frac{H_u}{\nu}+\nonumber\\
   &\sum_{t=1}^{K-2}\frac{\sqrt{(x_{u}[t+1]-x_{u}[t])^2+(y_{u}[t+1]-y_{u}[t])^2}}{\nu},
   \label{flighttime}
\end{align}

where $\frac{H_u}{\nu}$ is the distance to reach an altitude $H_u$ from the ground. We assume that this distance is achievable within one time slot of length $\tau$\footnote{For example, for $\nu=20m/s$, $H_u=50m$ and $\tau=4s$, $\frac{H_u}{\nu}\leq \tau$.}.

Our aim is to optimize UAVs' trajectories and  the subsets of visited IoT devices so that the AoU over devices is minimized. In the following, we first formulate our problem mathematically, and then we propose an approach to solve it efficiently.
% 
%Due to the limited wireless resources, we assume that only $N$ devices can transmit at a time. As a consequence, during a time interval $t$, a UAV can collect data from at most $N$ devices at a time.

%=============================================
% # III. Problem formulation #
% =============================================
\section{Problem Formulation }\label{Prob}
The main objective of this work is to minimize the AoU over devices during time span $T$. To this end, we optimize the stopping locations of UAVs over time along with the selected devices to collect data updates. Our optimization problem is formulated as follows.
 
\begin{mini!}%\label{Problem9}
{\boldsymbol{\alpha},\boldsymbol{x},\boldsymbol{y}}{\sum\limits_{i \in \mathcal{I}} f_i(\boldsymbol{\alpha}_i)
,\label{objective}}
{\label{GeneralOptimizati}}{}
% \addConstraint{ 
% R_{iu}[t]\geq \alpha_{iu}[t]R^{\rm min},\quad  \forall u \in \mathcal{U}, \forall i \in \mathcal{I}, \forall t \in \{0,\dots, K-1\}.
% \label{Rate}}
% {}{}
\addConstraint{R_{iu}[t]\geq \alpha_{iu}[t]R^{\rm min},\;  \forall u \in \mathcal{U}, \forall i \in \mathcal{I}, \label{Rate}}{}
\addConstraint{}{  \forall t \in \{0,\dots, K-1\} }\nonumber
{}{}
\addConstraint{ 
\sum \limits_{u\in \mathcal{U}}\alpha_{iu}[t]\leq 1, \; \forall i \in \mathcal{I},  \forall t \in \{0,\dots, K-1\}
\label{Association2}}
{}{}
\addConstraint{ 
\zeta_u(\boldsymbol{x}_u,\boldsymbol{y}_u)\leq \zeta_u^{\rm max}, \;  \forall u \in \mathcal{U}
\label{Time}}
{}{}
\addConstraint{ 
0\leq x_u[t]\leq x^{\rm max}, \;  \forall u \in \mathcal{U}, t \in \{0,K-1\}
\label{xposition}}
{}{}
\addConstraint{ 
0\leq y_u[t]\leq y^{\rm max}, \;  \forall u \in \mathcal{U}, t \in \{0,K-1\}
\label{yposition}}
{}{}
\addConstraint{ 
\alpha_{iu}[t]\in \{0,1\}, \;  \forall i \in \mathcal{I}, \forall u \in \mathcal{U}, t \in \{0,K-1\}
\label{alpha}}
% {}{}
% \addConstraint{ 
% (x_u[K-1],y_u[K-1])=(x_u^{\rm CDS},y_u^{\rm CDS}), \quad  \forall u \in \mathcal{U}.
% \label{CDS}}
{}{}
\addConstraint{(x_u[K-1],y_u[K-1])=(x_u^{\rm CDS},y_u^{\rm CDS}) \label{CDS}}{}
\addConstraint{}{ \forall u \in \mathcal{U}. }\nonumber
{}{}
\end{mini!}

Constraint~(\ref{Rate}) ensures that the rate between any UAV and its served IoT devices is above a pre-defined threshold $R^{\min}$. Constraint~(\ref{Association2}) prevents devices to transmit data to more than one UAV at a time. Constraint~(\ref{Time}) ensures that each UAV $u$ does not exceed a maximum flight time, $\zeta_u^{\rm max}$, that corresponds to its energy budget. Constraints~(\ref{xposition}) and (\ref{yposition}) limit UAVs' movements to the studied area. Finally, constraint~(\ref{alpha}) designates the binary variables, while constraint~(\ref{CDS}) ensures that each UAV returns to its CDS at the end of its flight. 

The underlying optimization combines binary and continuous variables. It also involves nonlinear objective function and constraints. As a consequence, the studied problem is 
MINLP and is challenging to solve. We note that traditional MINLP optimization algorithms (e.g., Branch and Bound) are not suitable for this problem as they require an exponential time of convergence and assume a global knowledge of the environment and its dynamics. Furthermore, classical machine learning algorithms, in particular standard RL, fail to scale with the large dimension of the set of variables and are not adapted to multi-agent environments. In the following, we leverage the MAPPO framework to solve the studied problem. %We also propose an adapted distributed meta-training algorithm to train the UAVs to optimize AoU in unseen environments.  
% ==================
% # IV. Multi Agent Reinforcement Learning
% ==================
% \section{Proposed Approach}\label{algo}
\section{Multi-Agent Proximal Policy Optimization Algorithm for AoU optimization (MAPPO-AoU)}\label{algo}

To obtain the optimal trajectories and subsets of served IoT devices in a multi-UAV setup is a complex problem, mainly because the UAVs optimize a common objective (i.e., the global AoU), and their policies are interdependent. Specifically, a UAV decides to visit and collect data from a given device based on whether the data of the target device has been collected in the past (by another UAV) or not. This is particularly challenging since, in practice, it is difficult to implement a fully synchronized system where each UAV has a global knowledge about the entire environment including other UAVs' positions and their historical actions (devices they have visited in the past, the time of data collection, etc.). We  address this problem from a multi-agent learning perspective. In particular, we adopt MAPPO approach where UAVs are assumed to use local information only to take their decisions during the execution time whereas they use a centralized learning during the training.

%As a consequence, it is more suited in a real world application to assume a partial observation of the

%For instance, suppose UAV $1$ collected data from device $1$ at $t_1$, the other UAVs 
%In fact, UAVs have limited information about devices, their data update generation, and other UAVs trajectories. The question that arises here is 

%The reward gleaned by each agent in a multi-agent system is not only related to its own actions but also to those of other agents.  The selection of the optimal policy of the other agents will be affected by the change of the policy of an agent and the estimation of the value function will be imprecise, so it won't be easy to ensure the convergence of the algorithm. 
\begin{figure}[ht]
\centering
   \includegraphics[width=0.75\linewidth]{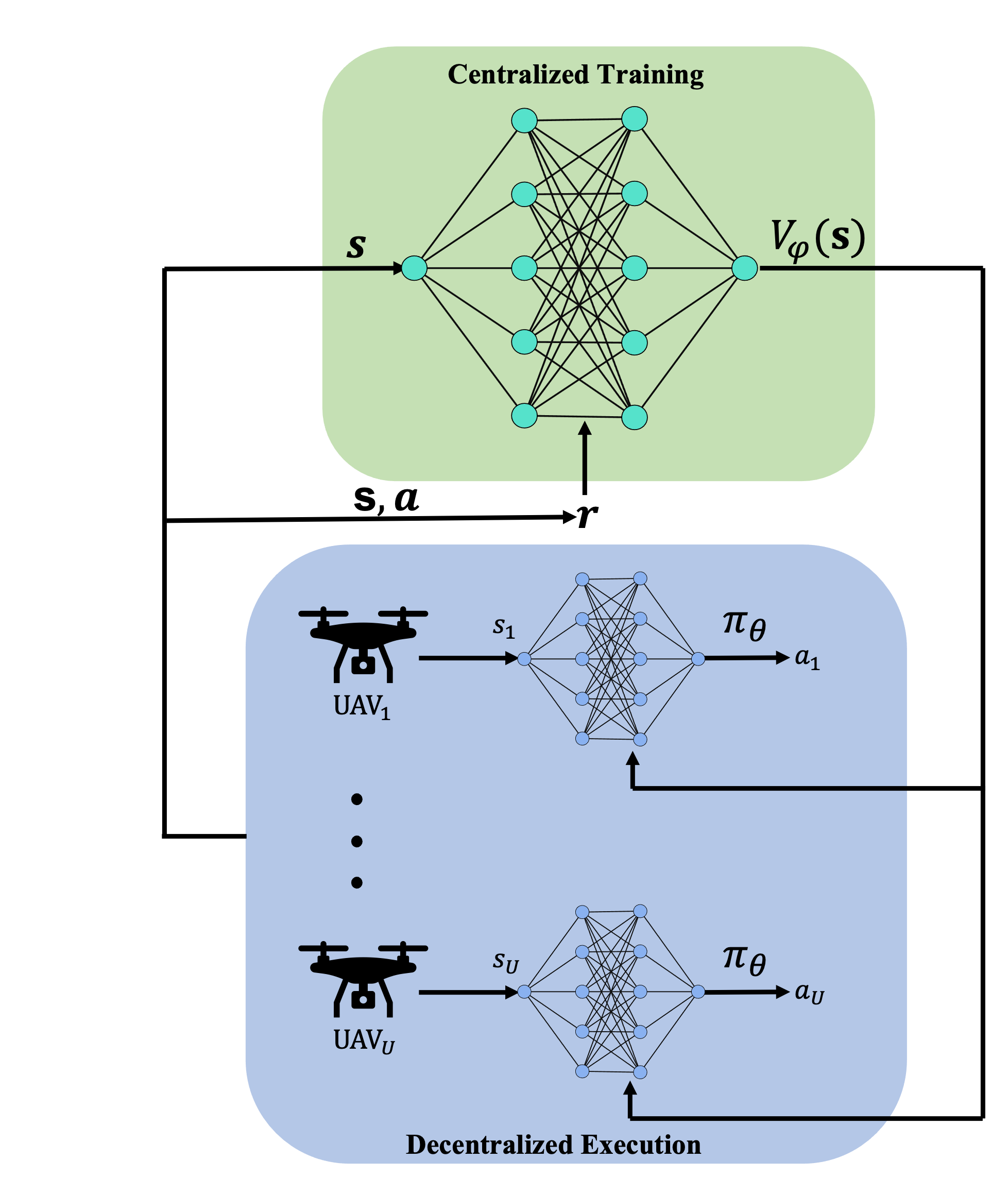}
   \caption{Centralized critic and decentralized actor for the training of MAPPO-AoU approach.}
    \label{fig:ctde}
\end{figure}

\begin{figure*}[h]
\minipage{0.3\textwidth}
   \includegraphics[width=1\linewidth]{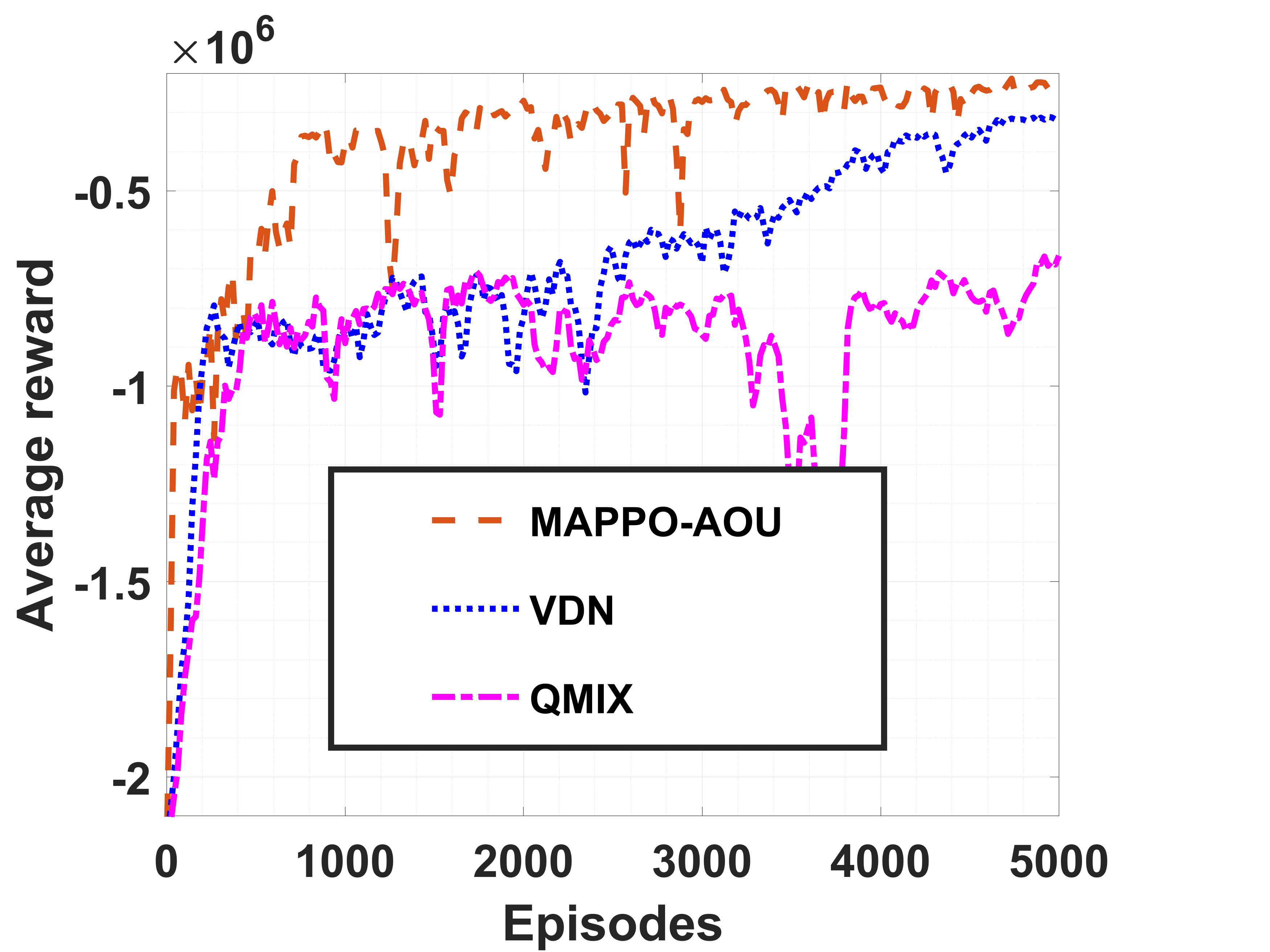}
   \caption{Average reward over episodes.}
    \label{fig:reward}
\endminipage\quad
\minipage{0.3\textwidth}
  \includegraphics[width=1\linewidth]{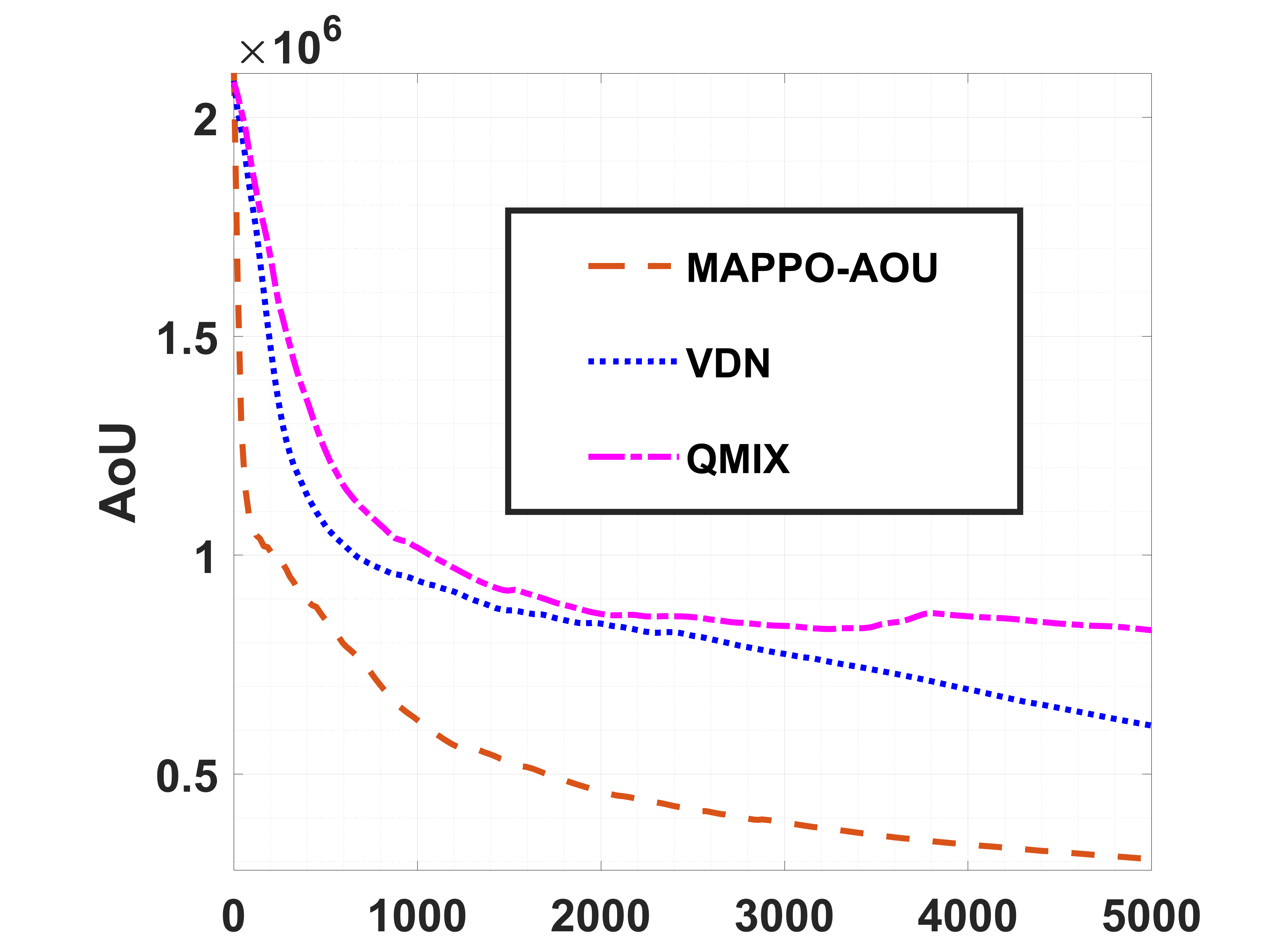}
\caption{AoU over episodes.}
 \label{fig:AoU}
\endminipage\qquad
\minipage{0.3\textwidth}
  \includegraphics[width=1\linewidth]{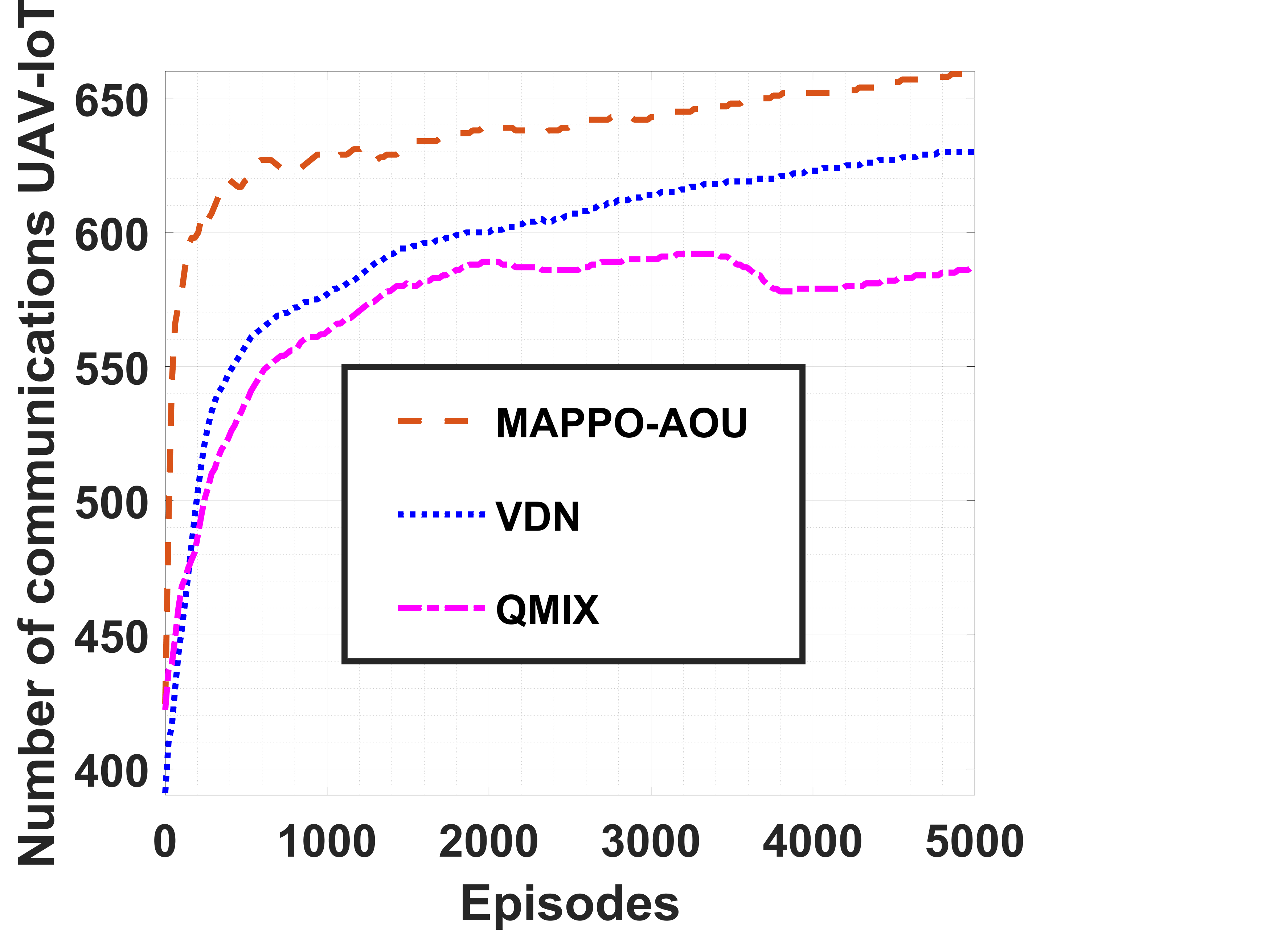}
   \caption{ Total number of UAV-IoT communications over time.}
    \label{fig:iotnondistinct}
    \endminipage
\end{figure*}

\begin{figure*}[ht]
\minipage{0.3\textwidth}
  \includegraphics[width=0.9\linewidth]{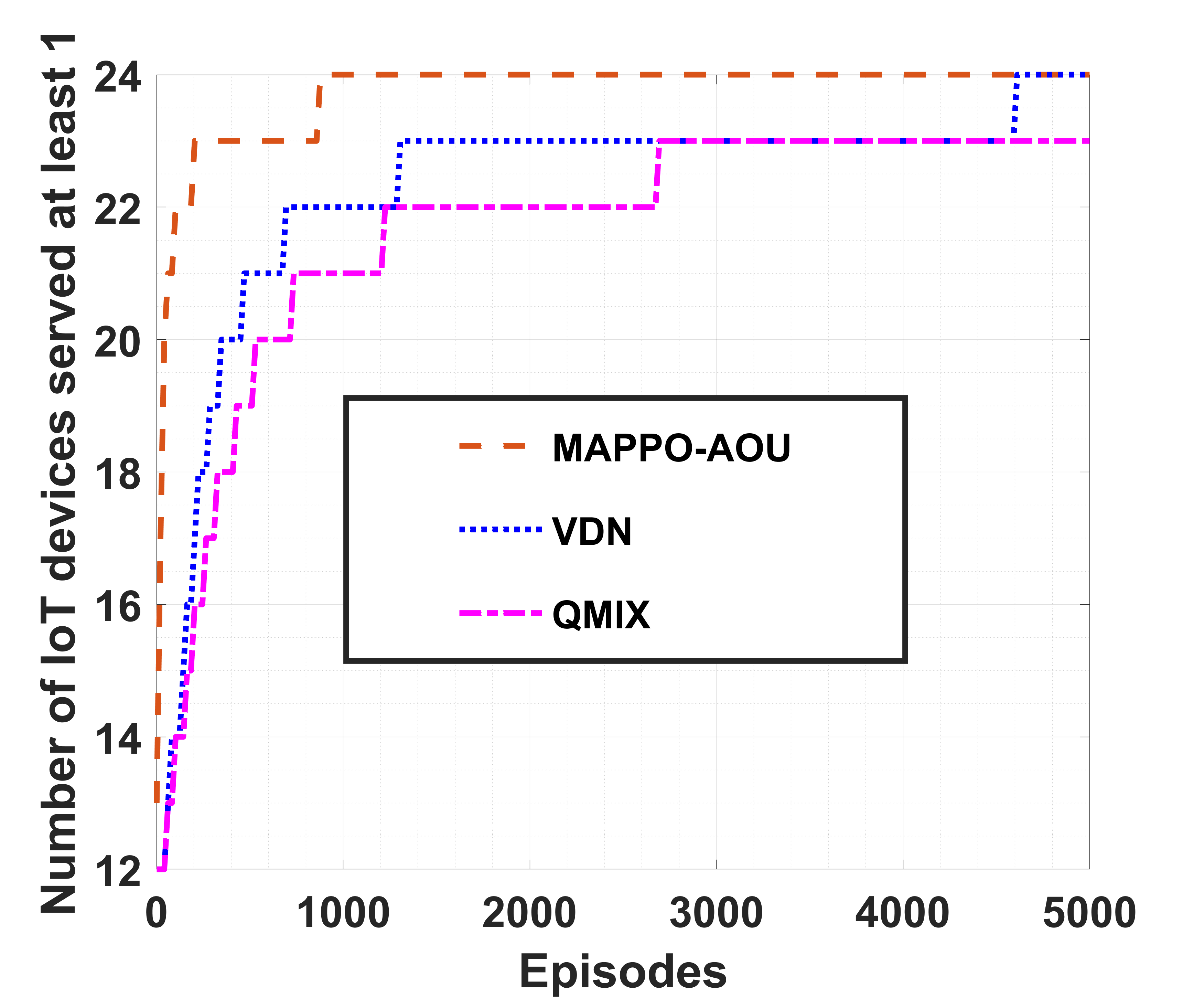}
   \caption{Number of served IoT devices at least \\ once over episodes.}
    \label{fig:iotdistinct}
\endminipage\qquad
\minipage{0.3\textwidth}
  \includegraphics[width=1\linewidth]{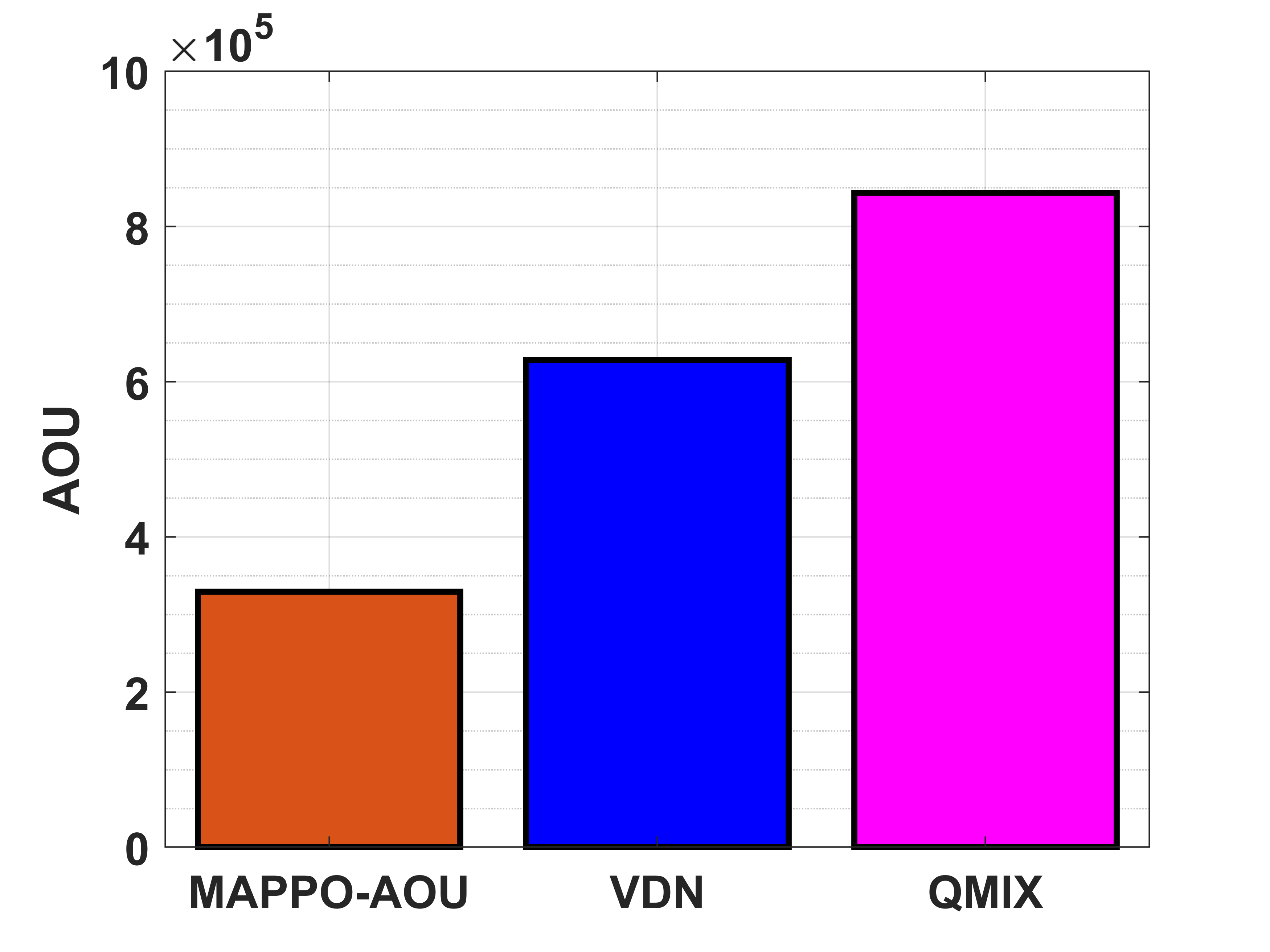}
   \caption{AoU value during the execution phase}
    \label{fig:testAoU}
\endminipage\qquad
\minipage{0.3\textwidth}
  \includegraphics[width=1\linewidth]{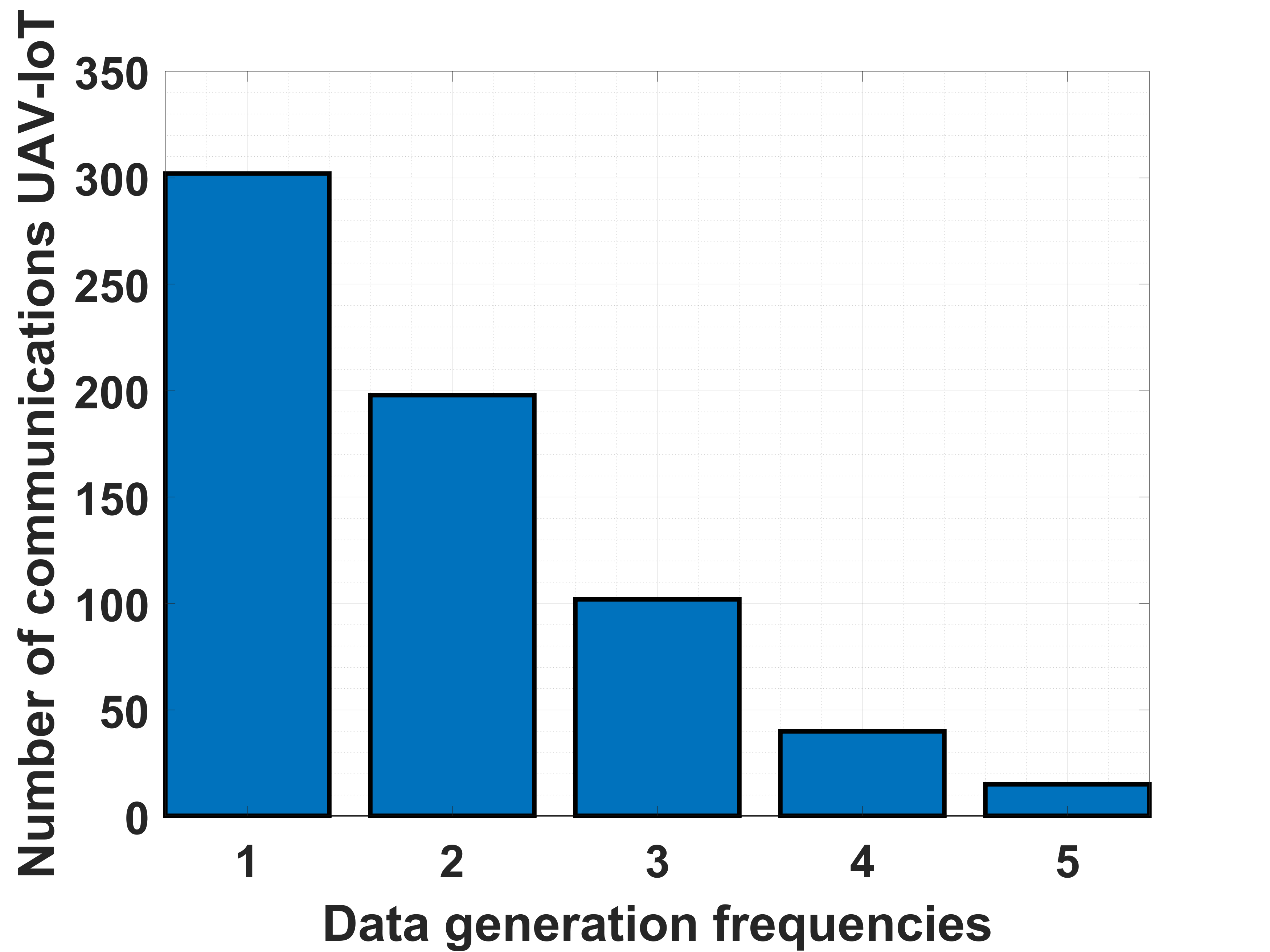}
   \caption{Number of served IoT devices over data \\ generation frequencies.}
    \label{fig:gendata}
    \endminipage
\end{figure*}
In order to solve this problem, MAPPO\cite{yu2021surprising} uses the idea of CTDE and extends the PPO  for the multi-agent system. In particular, it involves the actor-critic architecture that contains two components:
the actor that is responsible for taking decisions, and
 the critic evaluates the decisions by using a value
function. In a CTDE setup, the critic network learns a centralized value function that involves the knowledge of all agents' actions, while the actor network is used by each agent to calculate its policy only from its own local observations.

The proposed MAPPO-AoU algorithms consist of:
\begin{itemize}
    \item \textbf{Agents}: the set of UAVs.

    \item \textbf{Actions}: Each UAV has two decisions to take. First, it moves along the $x$ or $y$ axis by going \textit{up, down, right, left} or \textit{staying in the same place}. Second, it selects the subset of devices 
    from which it collects data. Thus, the action of a UAV $u$ at time $t$, $\boldsymbol{a}_u[t] \in \{0,1,2,3,4\}\times \{0,1\}^I$, where $0$ in $\{0,1,2,3,4\}$ denotes the action to stay at the same position, whereas $1,2,3,4$ designates respectively going \textit{up, down, right} or \textit{left}.
    
   % The action of each UAV is the UE that it seeks to collect data or return to CDS after completing its mission. Each UAV can move along the $x$ or $y$ axis. At time step $t$ the action taken by UAV $u$ is given by $a_u[t] \in \left\{ 0,1,2,3,4 \right\}$, where action $0$ is to stay as its position and action $1,2,3,4$ design respectively going  UP, DOWN, RIGHT, and LEFT. The joint action take by UAVs at timestep $t$ is, $\boldsymbol{a}[t]=[a_{u}[t], u=1, \ldots, U]$.
    
    \item \textbf{States}: The state of a UAV $u$ at time $t$ is given by its location and is denoted $s_u[t]$. To describe the set of all possible states, we discretize the studied area into a 2D grid with granularity $\delta$ that we denote $\mathcal{G}_\delta$. The set of states of a UAV is $\mathcal{S}=\mathcal{G}_\delta$. We use $\boldsymbol{s}[t]=(s_1[t],\dots,s_U[t])$ to denote the global state vector of all UAVs.
    %is the state of UAV $u$. Therefore, the set of states of all UAVs is $\mathcal{S}=\{\boldsymbol{s}[t], t=1\ldots T\}$, where $\boldsymbol{s}[t]=(s_1[t],\dots,s_U[t])$ is the global state vector. %The state of the UAV $u$ at the time step $t$ is $s_{u}[t]=(x_{u}[t],y_{u}[t])$ and the global state of UAVs at timestep $t$ is $\boldsymbol{s}[t]=[{s}_{u}[t], u=1, \ldots ,U]$.
    
 %   At time step $t$ the next state is computed as follows, \begin{equation}
%s_u[t+1]= \begin{cases}s_u[t]+\delta, & \text { if } a_{u}[t]=1, \\ s_u[t]-\delta, & \text { if } a_{u}[t]=2, \\ s_u[t]+\delta, & \text { if } a_{u}[t]=3, \\ s_u[t]-\delta, & \text { if } a_{u}[t]=4, \\ s_u[t], & \text { otherwise. }\end{cases}
%\label{UAV_dynamic}
%\end{equation}
If a UAV reaches the edge of the grid at time step $t$ and $a_u[t]$ would cause it to leave the grid in $t+1$, it will stay in its current location at $t+1$.
    
    \item  \textbf{Rewards}: The reward is closely related to the total AoU of all UAVs as defined in (\ref{objective}). The reward at time step $t$ is given by
    \begin{equation}
        r(\boldsymbol{a}[t], \boldsymbol{s}[t])=-\sum\limits_{i=0}^{I}\sum\limits_{n=0}^{\lfloor \frac{K}{k_i} \rfloor}A_{i}^n[t]
        \label{indiv_reward}
    \end{equation}
    
    \item \textbf{Policy function}: We define a policy function $\boldsymbol{\pi}_{\boldsymbol{\theta}}(a_{u},s_{u})$ represented by the actor network and parameterized by the vector $\boldsymbol{\theta}$, it is used to generate the UAV strategy based on its own local state.% The agent shares the same network but each UAV $u$ takes action based on its own local state $s_{u}$, and optimizes its accumulated discounted accumulated reward.
    % $J(\theta)=\mathbb{E}_{a[t], s[t]}\left[\sum_t \gamma^t R\left(s[t], a[t]\right)\right]$.

    \item \textbf{Value function}: The centralized value function $V_{\boldsymbol{\phi}}(\boldsymbol{s})$ is represented by the critic network and parameterized by the vector $\boldsymbol{\phi}$. This network is used to estimate the achievable future reward by UAVs at every state $\boldsymbol{s}$. The goal of the critic network is to find the global optimal strategy that maximizes the expected reward. 
\end{itemize}
We note that $V_\phi(\boldsymbol{s})$ is used in the training only. As a consequence, it can take global information as input. 

As illustrated in Fig.~\ref{fig:ctde}, during the training, MAPPO-AoU alternates between the actor and critic until it achieves convergence. Specifically, UAVs locally train the policy actor using a PPO approach. At each training step $t$, and given a UAV state $s_u[t]$, the actor network is trained, and action is sampled according to the policy $\pi_\theta(a_u[t],s_u[t])$. Then the reward can be observed by executing the joint action based on the current policy, subject to satisfying the constraints (\ref{Association2}) and (\ref{Rate}). Finally, the resulted global state vector of all UAVs is then fed to the critic network which in turn is trained by minimizing a loss function and using a PPO approach. We mention that PPO is a policy gradient method that uses a clipped surrogate objective and handles large ML models efficiently. Due to space limitations, the details about PPO are removed. For reference, the reader can check the work in~\cite{schulman2017proximal}. 

\section{Simulation Results}\label{simu}
\begin{table}[ht]
\centering
\begin{tabular}{|l|l||l|l|}
\hline
\textbf{Parameter} & \textbf{Value} & \textbf{Parameter} & \textbf{Value}   \\ \hline
$I$       & $25$      &$U$       & $3$  \\ \hline
$x_{max}$       & $1000m$ & $y_{max}$       & $1000m$  \\ \hline
$H_u$       & $[80,100]m$   & $R_{min}$       & $150$ Kbit/s  \\ \hline
$\tau$       & $3ms$   &$k_i$       & $[1,5]$  \\ \hline
$B_{i,u}$       & $[1.5,2]$ GHz   &$\sigma^{2}$       & $-120dBm$  \\ \hline
$P_{i}$       & $[0,1]$ mW   &$T$       & $500\tau$  \\ \hline
%$\nu$       & $15m/s$ \\ \hline
\end{tabular}
\caption{Experiment setup}
\label{tab:my-table-parameter}
\end{table}
% \begin{table}[ht]
% \centering
% \begin{tabular}{|l||l|}
% \hline
% \textbf{Parameter} & \textbf{Value}   \\ \hline
% $I$       & $25$  \\ \hline
% $U$       & $3$  \\ \hline
% $x_max$       & $1000m$  \\ \hline
% $y_max$       & $1000m$  \\ \hline
% $H_u$       & $[80m,100m]$  \\ \hline
% $R_{min}$       & $150000$  \\ \hline
% $\tau$       & $3ms$  \\ \hline
% $k_i$       & $[1,5]$  \\ \hline
% $\beta_{0}$       & $1$  \\ \hline
% $\alpha$       & $2$  \\ \hline
% $B_{i,u}$       & $[1.5MHz,2MHz]$  \\ \hline
% $\sigma^{2}$       & $-120dBm$  \\ \hline
% $P_{i}$       & $[0,1]$  \\ \hline
% $T$       & $500\tau$  \\ \hline
% $V$       & $15m/s$  \\ \hline
% \end{tabular}
% \caption{Experiment setup}
% \label{tab:my-table-parameter}
% \end{table}

To evaluate the performance of the proposed approach, we  consider an area of $1000m\times1000m$, where a number of $25$ IoT devices are randomly scattered. We also suppose that $3$ UAVs are deployed to collect the data and keep it as fresh as possible. The UAVs hover at altitudes between $80m$ and $100m$, with a constant speed equal to $15m/s$. The mission time is assumed to be equal to $\zeta_u^{\rm max}=500\tau$, where $\tau=3ms$. Moreover, the data generation frequency of IoT devices is randomly chosen between $1\tau$ and $5\tau$. Devices are assigned a fixed bandwidth, randomly picked between $[1.5,2]$ GHz and a constant power between $[0,1]$ mWatt. To satisfy the quality of service constraint, the minimum rate is set to $150$ Kbit/s. The parameters of our simulation setup are summarized in Table \ref{tab:my-table-parameter}.

As benchmarks, we use the popular off-policy MARL algorithms; VDN ~\cite{sunehag2017value} and QMIX algorithms\cite{rashid2018qmix}. VDN is a framework for cooperative MARL tasks in which each agent learns a separate value function that is decomposed into shared and local value functions. The shared value function is then used to train the agents to work together, whereas the local value function allows each agent to learn its own policy. QMIX is a deep multi-agent reinforcement learning algorithm that enables cooperative learning by using a monotonic factorization of the value function. %It has been used to solve a variety of tasks and is based also on the CTDE paradigm. QMIX employs a mixture network that combines individual agent values into a single global value. This value is then used to determine the best action to take for each agent. 

%As benchmarks, we use two off-policy MARL algorithms, the VD Network (VDN) algorithm\cite{sunehag2017value} and the QMIX\cite{rashid2018qmix} algorithm. VDN is a framework for cooperative MARL tasks that follows the CTDE model. It is a value-based learning approach in which each agent learns a separate value function that is decomposed into a shared and a local value function. The shared value function is then used to train the agents to work together, whereas the local value function allows each agent to learn its own policy. VDN assumes that the total reward has an additivity property and that the agents' policies, given the joint state, should be conditionally independent. QMIX is a deep multi-agent reinforcement learning algorithm that enables cooperative learning by using a monotonic factorization of the value function. It has been used to solve a variety of tasks and is based also on the CTDE paradigm. QMIX employs a mixture network that combines individual agent values into a single global value. This value is then used to determine the best action to take for each agent. 

The results presented in \Cref{fig:reward,fig:AoU,fig:iotdistinct,fig:iotnondistinct,fig:testAoU,fig:gendata} demonstrate the superiority of MAPPO-AoU approach over the benchmark algorithms. In Fig.~\ref{fig:reward}, we plot the average reward gleaned by UAVs over episodes. The plot shows that MAPPO-AoU converges faster and achieves higher rewards than the benchmarks. Fig. \ref{fig:AoU} depicts the AoU of the system over episodes during the training. It shows that the MAPPO-AoU approach results in the UAVs discovering optimal strategies to minimize the AoU. The results of this figure are in line with some recent work that shows (empirically) that PPO achieves high performance in cooperative games compared to off-policy approaches~\cite{Yu2021TheSE}.%\cite{Yu2021TheSE, DBLP:journals/corr/abs-2006-07869}.

In Fig. \ref{fig:iotnondistinct}, we plot the total number of communications between IoT devices and UAVs during the time mission. This number represents the total number of times an IoT device sent its collected data to a UAV. 
Moreover, in Fig. \ref{fig:iotdistinct}
, we plot the number of IoT devices that have been visited at least once. These two figures highlight that the MAPPO-AoU solution empowers UAVs to find the most effective strategy for visiting more IoT devices and reducing the AoU. In fact, to reduce the total AoU of the network, UAVs are required to visit a maximum number of devices and collect their data. This is because an unvisited device for a long time can increase the AoU excessively, resulting in a sub-optimal solution.

In Fig. \ref{fig:testAoU}, we conduct an experiment using our trained agents to measure the AoU during the execution. As one may expect, Fig. \ref{fig:testAoU} indicates that the minimum value of AoU is achieved through MAPPO-AoU. Compared to VDN and QMIX, the AoU during the execution is reduced by a factor of $1/2$ and $2/3$ respectively. Finally, in Fig. \ref{fig:gendata},
we plot the number of UAV-IoT communications against the data generation frequencies. As shown in the figure, devices with high generation frequencies result in a higher number of communications with UAVs. This is simply due to the fact that devices with high data generation rates are visited more often compared to devices with lower data generation.

%trained UAVs using MAPPO-AoU have discovered their optimal policy for visiting more IoT devices with low data generation frequencies, leading to better minimization of AoU.

% The analysis of these results revealed  the effectiveness of our proposed approach. in fact, these results show that our proposed solution MAPPO-AOU converges faster and achieves higher rewards compared to the benchmark algorithm. Additionally, the results demonstrate that our MAPPO-AOU is able to quickly find a global strategy that minimizes the AoU. Also using our approach UAVs are able to find the best strategy to visit more IoT devices. Finally, our proposed approach enables the UAVs to discover the optimal policy that allows them to serve IoT devices with lower data generation frequency. 

% ==================
% V. Conclusion #
% ==================

\section{Conclusion}\label{Conc}

In this work, we studied the deployment of UAVs as data collectors. In particular, we investigated the case where the collected data is time-sensitive. To measure the data freshness, we introduced the AoU metric and modeled the studied problem as an MINLP. To efficiently solve our optimization problem, we proposed a policy based MARL approach. Our simulation results show that the proposed approach outperforms conventional off policy based RL methods. 

%In this work, we present a comprehensive framework for cooperative UAVs to maximize the freshness of collected data. We measured data freshness with a metric called AoU and enhanced it by including the data generation frequency of IOT devices. The problem is formulated as a non-convex optimization problem, and we propose a MARL-based solution to solve it. Our proposed solution achieves faster convergence and better performance than the benchmark algorithm, according to simulation results. 

% Additional research, such as optimizing the 3D positions of UAVs by adapting the algorithm to continuous action spaces can be done.

% ==============
% # REFERENCES #
% ==============
\balance

\bibliographystyle{IEEEtran}
\bibliography{IEEEabrv,biblio_traps_dynamics}

\end{document}